# The Riley slice revisited

Yohei Komori
Caroline Series

**Abstract** In [4], Keen and Series analysed the theory of pleating coordinates in the context of the Riley slice of Schottky space $\mathcal{R}$, the deformation space of a genus two handlebody generated by two parabolics. This theory aims to give a complete description of the deformation space of a holomorphic family of Kleinian groups in terms of the bending lamination of the convex hull boundary of the associated three manifold. In this note, we review the present status of the theory and discuss more carefully than in [4] the enumeration of the possible bending laminations for $\mathcal{R}$, complicated in this case by the fact that the associated three manifold has compressible boundary. We correct two complementary errors in [4], which arose from subtleties of the enumeration, in particular showing that, contrary to the assertion made in [4], the *pleating rays*, namely the loci in $\mathcal{R}$ in which the projective measure class of the bending lamination is fixed, have two connected components.

In [4], L Keen and C Series used their theory of pleating invariants to study the so called *Riley slice of Schottky space*. The Riley slice is shown to be foliated by *pleating rays* on which the geometry of the limit set has fixed combinatorial properties. There are two rather subtle errors in [4], concerning the labelling and connectivity of these rays. While the errors do not substantially affect the main results, the correct picture illustrates the interesting new phenomenon that pleating varieties may not be connected, as well as some delicate points related to the marking of the group. The ideas involved will apply to other examples, and we consider them to be of sufficient interest to be worth discussing at some length. Besides explaining and correcting the errors, we take the opportunity to review the background and discuss some of the techniques used in [4] in more detail.

A Kleinian group is a discrete subgroup $G$ of $PSL(2, \mathbf{C})$. It acts on the Riemann sphere $\hat{\mathbf{C}}$ by Möbius transformations and on hyperbolic 3–space $\mathbf{H}^3$ by isometries. The regular set $\Omega(G)$ is the subset of $\hat{\mathbf{C}}$ on which the elements of $G$ form a normal family and the limit set $\Lambda(G)$ is its complement. The quotient $\Omega(G)/G$ is a (possibly disconnected) Riemann surface and $\mathbf{H}^3/G$ is a hyperbolic 3–manifold whose ends are exactly the components of $\Omega(G)/G$. Let $\mathcal{C}$ be the hyperbolic convex hull of $\Lambda(G)$ in $\mathbf{H}^3$; $\mathcal{C}/G$ is the *convex core* of the hyperbolic manifold $\mathbf{H}^3/G$. The boundary $\partial \mathcal{C}/G$ of $\mathcal{C}$ the convex core of $\mathbf{H}^3/G$ is a (possibly disconnected) pleated surface homeomorphic to $\Omega(G)/G$. We denote the geodesic lamination along which this surface is pleated by $pl(G)$.

Let $G_\mu$ be a family of Kleinian groups depending holomorphically on a parameter $\mu$ which varies over a complex manifold $D$, and such that the groups $G_\mu$ are all quasiconformally conjugate. The theory of pleating invariants analyses $D$ in terms of the *pleating varieties* $\mathcal{P}_\lambda = \{\mu \in D : pl(G_\mu) = \lambda\}$, where $\lambda$ is a fixed geodesic lamination on $\partial \mathcal{C}/G$.

Let $f\colon U \mapsto \mathbf{C}$ be a holomorphic function defined on a subset $U \subset D$. The *real locus* of $f$ in $U$ is the set $f^{-1}(\mathbf{R}) \cap U$. A geodesic lamination is called *rational* if all its leaves are closed.

In all cases studied so far, [3, 4, 5], $D$ has one or two complex dimensions and it has been shown that:

(1) All geometrically possible pleating varieties are non-empty.

(2) The pleating variety $\mathcal{P}_\lambda$ is a union of connected components of the real loci of a (finite) collection of non-constant holomophic functions $f_{i,\lambda}$ in the part of $D$ on which $G_\mu$ is non-Fuchsian.

(3) The pleating variety $\mathcal{P}_\lambda$ is a submanifold of appropriate dimension.

(4) The pleating varieties $\mathcal{P}_\lambda$ for which $\lambda$ is rational are dense in $D$.

The pleating varieties foliate $D$, possibly omitting an exceptional set on which $G_\mu$ is Fuchsian.

We say that $g_1, g_2 \in G$ are *I–equivalent* if $g_1$ is conjugate in $G$ to either $g_2$ or $g_2^{-1}$, and write $g_1 \sim g_2$. We denote the equivalence class of $g$ by $C(g)$ and note that the trace function $\operatorname{Tr} g$ is constant on $C(g)$. An oriented closed geodesic in $\mathbf{H}^3/G$ corresponds to a conjugacy class in $G$, however a closed leaf of a geodesic lamination is unoriented and hence defines only an an I–equivalence class in $G$.

Suppose that the lamination $\lambda$ is rational. The functions $f_{i,\lambda}$ of (2) above may be taken to be the set of trace functions $\operatorname{Tr} g_{i,\lambda}$ as $g_i = g_{i,\lambda}, i = 1, \ldots, k$ ranges over a full set of representatives of the I–equivalence classes corresponding to leaves of $\lambda$. All of these trace functions are, in principle, computable holomophic functions on the parameter space $D$. (In any specific example we have to discuss how to make a consistent choice of sign for the trace corresponding to a lifting from $PSL(2, \mathbf{C})$ to $SL(2, \mathbf{C})$; in the case of this paper the problem does not arise since the deformation space is defined as a set of subgroups of $SL(2, \mathbf{C})$.) It follows that, in order to find the foliation by pleating varieties, and hence to compute $D$, it suffices to enumerate the possible rational pleating laminations $\mathcal{P}_\lambda$ for $\partial \mathcal{C}/G$, and then to identify $\mathcal{P}_\lambda$ among the components of the real loci of the associated trace functions $\operatorname{Tr} g_{i,\lambda}$.

In the case of the Riley slice, this programme was carried out in [4].

Consider the set of discrete subgroups of $SL(2, \mathbf{C})$ which are freely generated by two non-commuting parabolics. Up to conjugation in $SL(2, \mathbf{C})$, any such subgroup can be put in the form $G = G_\rho = \langle X, Y_\rho \rangle$ where $X, Y_\rho$ are the matrices
$$X = \begin{bmatrix} 1 & 1 \\ 0 & 1 \end{bmatrix}, \ Y_\rho = \begin{bmatrix} 1 & 0 \\ \rho & 1 \end{bmatrix}.$$

The Riley slice $\mathcal{R}$ is defined by:

$$\mathcal{R} = \{\rho \in \mathbf{C} : \Omega(G_\rho)/G_\rho \text{ is a four times punctured sphere }\}.$$

It is known that the deformation space $D = \mathcal{R}$ is topologically an annulus in $\mathbf{C}$ (see [2]). Thus the real locus of a holomorphic function on $\mathcal{R}$ has one real dimension. In this case, (2) is theorem 3.7 in [4], (1) and (3) follow from theorem 4.1 and (4) is theorem 5.2.

The first error in [4] concerns the enumeration of rational laminations on $\partial \mathcal{C}/G$. Let $S$ denote a four times punctured sphere. In the Riley slice setup, the convex hull boundary $\partial \mathcal{C}/G$ is compressible, in other words, the induced map $\pi_1(S) \to G$ is not injective. This means that to enumerate correctly the possible pleating laminations, one has to determine when two distinct laminations on $S$ define the same family of geodesics in $\mathbf{H}^3/G$. In particular, one has to determine when the images of distinct simple closed curves on $S$ are equal in $G$. This was not handled quite correctly in [4]. Curves counted there as distinct are actually equivalent in pairs. (For a more general discussion of the pleating locus with compressible boundary, we refer to Otal's thesis [9].)

The second error concerns connectivity of the pleating varieties; contrary to the assertion in [4], except for degenerate cases, each pleating variety has *two* connected components. In fact, the inductive proof of theorem 4.1 of [4] is not quite correct.

At the time of writing [4], these errors went unnoticed as they are in some sense complementary. That errors were present, was first deduced as follows. Let $\gamma$ be a simple closed curve on $S$ and let $\rho^*(\gamma)$ be the endpoint on $\partial \mathcal{R}$ of the pleating ray $\mathcal{P}_\gamma$, which we assume for the moment has one connected component as in theorem 4.1 of [4]. In the group $G_{\rho^*(\gamma)}$, the element $\gamma$ is represented by an accidental parabolic with trace equal $-2$ (see [4] proposition 4.2). Since the trace is a polynomial in $\rho$ with integer coefficients, the element $\gamma$ in the group $G_{\overline{\rho^*(\gamma)}}$ is also an accidental parabolic so that one would expect $\overline{\rho^*(\gamma)}$ also to be an endpoint of $\mathcal{P}_\gamma$. This appears to contradicts the assertion of theorem 4.1 that $\mathcal{P}_\gamma$ has a unique branch. Again, according to [4], the possible rational pleating laminations are enumerated by the rationals modulo 2 and the points $\rho^*(\gamma), \overline{\rho^*(\gamma)}$ should be the endpoints of a pair of *distinct* pleating laminations whose labels are $p/q$ and $2 - p/q$. However, this would imply that the distinct elements $\gamma(p/q)$ and $\gamma(2 - p/q)$ are both pinched at $\rho^*(\gamma)$ (and also $\overline{\rho^*(\gamma)}$), which is impossible.

These contradictions are resolved simultaneously by showing that in fact (a) the laminations with labels $p/q$ and $2-p/q$ are the same, and (b) the pleating locus has two connected components which are complex conjugate in the $\rho$–plane.

The details of how this works are explained below.

## 1   Enumeration

As explained above, to enumerate correctly the possible rational pleating laminations, one has to determine when two distinct homotopy classes of simple closed curves on $S$ define the same geodesic in $\mathbf{H}^3/G_\rho$. This involves an implicit choice of marking on $S$ (i.e. a choice of generators for $\pi_1(S)$), together with a choice of homomorphism $h\colon \pi_1(S) \to G$. To explain the error, we have to review the enumeration with some care.

The group $G_\rho = \langle G_\rho; X, Y_\rho \rangle$ should be thought of as *marked* by the ordered pair of generators $(X, Y_\rho)$. Thus, although $Y_\rho^{-1} = Y_{-\rho}$ so that $G_\rho = G_{-\rho}$ as subgroups of $PSL(2, \mathbf{C})$, the marked groups $\langle G; X, Y_\rho \rangle$ and $\langle G; X, Y_\rho^{-1} \rangle$ are distinct. (In fact, it follows easily from lemma 1 of [2], that the only possible pair of parabolic generators of $G_\rho$ are the (unordered) pair $X^{\pm 1}, Y_\rho^{\pm 1}$ .) Thus we should always identify $G_\rho$ and $G_{\rho'}$ by the isomorphism $X \mapsto X, Y_\rho \mapsto Y_{\rho'}$. We denote by $\langle G; X, Y \rangle$ an abstract two generator marked free group and always use the isomorphisms $X \mapsto X, Y \mapsto Y_\rho$ to identify $G$ with $G_\rho$. With these identifications understood, abstract words in the symbols $X^\pm, Y^\pm$ represent elements in both the groups $G$ and $G_\rho$.

The content of proposition 2.1 of [4] is that a line of rational slope $p/q \in \mathbf{Q} \cup \infty$ in $\mathbf{C}$ projects to a homotopy class of simple closed non-boundary parallel curves on $S$, and that every such homotopy class on $S$ is obtained in this way. Further, the homotopy classes corresponding to distinct rationals are distinct.

Given a hyperbolic structure on $S$, there is a unique closed geodesic in each free homotopy class of simple closed non-boundary parallel curves. In [4], with the hyperbolic structure of $\partial \mathcal{C}/G$ understood, we denoted the geodesic corresponding to a line of slope $p/q$ by $\gamma(p/q)$. There is some confusion at this point, which contributes to the error under discussion. Since the line described in proposition 2.1 is in fact unoriented, the correct statement is that a line of rational slope defines an I–equivalence class in $\pi_1(S)$. From now on, therefore, $\gamma(p/q)$ should be understood to denote a specific I–equivalence class in $\pi_1(S)$.

The essence of the proof of proposition 2.1 appears on page 78 of [4] as part of the proof of proposition 2.2, which describes an explicit word $V_{p/q}$ in the generating set $\{X^{\pm 1}, Y^{\pm 1}\}$ which represents $h(\gamma(p/q))$ in $G$. We need to go through the construction of $V_{p/q}$ with some care. The idea is a simple case of the method of $\pi_1$–train tracks introduced in [1], see also [7].

Following [4], let $\mathcal{L}$ denote the integer lattice in the complex plane $\mathbf{C}$; let $\beta \colon z \mapsto z + 2i$ and let $\xi, \eta$ be the rotations by $\pi$ about the points $i$ and $i+1$ respectively. Let $\Gamma_0 = \{\xi^\pm, \eta^\pm, \beta^\pm\}$. The surface $S$ can be realised as the quotient of $\mathbf{C} - \mathcal{L}$ by the group $\Gamma$ generated by the elements of $\Gamma_0$.

We shall compare the three diagrams in figure 2 in [4]. Figure 2a is a fundamental domain $R$ for the action of $G = G_\rho$ on $\Omega = \Omega_\rho$. Figure 2b is a fundamental domain $R'$ for the action of $\tilde{G} = \tilde{G}_\rho$, the Fuchsian uniformisation of $\pi_1(S)$, acting in the hyperbolic disc $\Delta$, thought of as the universal cover of $S$. Figure 2c is a fundamental domain $R''$ for the action of $\Gamma$ on $\mathbf{C} - \mathcal{L}$. The sides of each of these domains are supposed to be labelled by generators $\alpha$ of $G$, $\tilde{G}$ and $\Gamma$ respectively in such a way that the label $\alpha$ on a side indicates that it is paired to the side labelled $\alpha^{-1}$ under the action of $\alpha$. (We note that, although $R''$ is a rectangle in $\mathbf{C}$, for the purposes of this discussion it should be thought of as having *six* sides.) Denote by $\tilde{G}_0$ the generating set $\{X', X'^{-1}, Y', Y'^{-1}, B', B'^{-1}\}$ of $\tilde{G}$.

Unfortunately, there is a labelling error in figure 2 (but not in the text) which may obscure the explanation on pages 77–78 of [4]. The configuration in figure 2(a) refers to the case $\rho < -4$. The two circles shown are the isometric circles of $Y_\rho^\pm$; since $\rho < -4$ the circle on the right has centre $-1/\rho$ and is the isometric circle of $Y_\rho$. This circle is identified with the circle on the left by $Y_\rho$ and thus, with the convention explained above, the labels $Y$ and $Y^{-1}$ should be interchanged. This error carries through to figures 2 (b) and (c) in which we should interchange the labels $Y'$ and $Y'^{-1}$, and $\eta$ and $\eta^{-1}$, respectively.

We proceed with this change of labelling throughout.

Let $\gamma$ be *any* simple closed non-boundary parallel loop on $S$. Its lift to any of the three covering spaces $\Omega$, $\Delta$ or $\mathbf{C}-\mathcal{L}$ of $S$ is simple and therefore appears on each region $R, R', R''$ as a collection of pairwise disjoint arcs with endpoints on the labelled sides. When the sides of one of the regions $R, R', R''$ are identified by the side pairings, there is a unique way to link the endpoints of the arcs to form a simple closed loop on $S$. This loop is well defined up to orientation and homotopy, and thus we obtain an I–equivalence class in $\pi_1(S)$.

Making a suitable homotopy, we may assume that none of these arcs join a side to itself. It is also clear that the total number of arcs meeting a side labelled $\alpha$ must equal the number meeting its paired side labelled $\alpha^{-1}$. Let $n(\alpha, \beta)$ denote the number of arcs joining the sides with labels $\alpha, \beta$. When the sides of $R'$ are identified, any arc joining sides $X'$ to $X'^{-1}$ links up to form a loop round a puncture. Since $\gamma$ is connected and non-boundary parallel, we conclude that $n(X', X'^{-1}) = 0$, and likewise that $n(Y', Y'^{-1}) = 0$. A similar argument (which makes crucial use of the fact that $\gamma$ is simple) shows that at least one of $n(X', B')$ and $n(X'^{-1}, B'^{-1})$, and at least one of $n(Y', B')$ and $n(Y'^{-1}, B'^{-1})$, must vanish, see [1, 7].

Exactly the same constraints apply to the weights $n(\alpha, \beta)$, $\alpha, \beta \in \Gamma_0$, in figure 2c. Inserting these constraints, we obtain precisely either one of the three patterns shown in figure 3 of [4], or its reflection in the line $\Re z = 1/2$. In these diagrams, there is at most one line $l(\alpha, \beta)$ joining a pair of sides $\alpha, \beta$ and the integer label $k$ on $l(\alpha, \beta)$ indicates that $n(\alpha, \beta) = k$. Conversely, given such a weighted diagram, we can recover a simple closed curve by replacing the line $l(\alpha, \beta)$ by $n(\alpha, \beta)$ parallel arcs joining the sides $\alpha, \beta$. When the sides of $R''$ are identified, there is a unique way to link these arcs to form a union of simple closed loops on $S$. There is one connected loop if and only if the integers $n, m$ appearing in figure 3 are relatively prime. Taking $(n, m) = 1$, we see that each of the patterns in figure 3 is exactly that obtained from a line of rational slope in $\mathbf{C}$, and that, provided we include reflections as above, every line of rational slope appears. This proves that every homotopy class of simple closed loops on $S$ is the projection of a line of rational slope in the plane as claimed.

We want to show that lines of different slope correspond to non-homotopic loops on $S$. To do this, observe that each side $\sigma$ of $R'$ is a line joining two punctures on $S$, and that the number of arcs meeting $\sigma$ is exactly the minimum geometric intersection number of loops homotopic to $\gamma$ with $\sigma$. It is clear that these intersection numbers determine $n(\alpha, \beta)$, $\alpha, \beta \in \tilde{G}$, which gives the result. It is also clear from the weighted diagrams in figure 3, that *all* lines of the same slope define the same class.

As noted above, this construction determines a curve only up to homotopy and orientation. The class corresponding to a line of slope $p/q$ in figure 2c is exactly the I–equivalence class $\gamma(p/q)$ in $\pi_1(S)$ described above.

We now want to find a word $V_{p/q}$ representing $h(\gamma(p/q))$ in $G_\rho$. As indicated in the proof of 2.2 in [4], this is done by the method of cutting sequences, see for example [1, 10]. We explain the method in somewhat more detail here.

Consider first the tesselation $\mathcal{T}$ of the hyperbolic disc $\Delta$ by images of the region $R'$ under the action of the group $\tilde{G}$. With the correction noted above, the sides of $R'$ should be labelled, in anticlockwise order starting from 0 by $B', Y'^{-1}, Y', B'^{-1}, X'^{-1}, X'$. These labels are transported to the tesselation $\mathcal{T}$ by the action of $\tilde{G}$. Two copies $R'_1, R'_2$ of $R'$ meet along each edge, and each edge carries two labels $\alpha, \alpha^{-1} \in \tilde{G}_0$, one label interior to $R'_1$ and the other interior to $R'_2$.

Let $\lambda$ be an oriented geodesic segment in $\Delta$ and let $\alpha_1, \ldots, \alpha_k$ be the ordered sequence of labels of edges of $\mathcal{T}$ cut by $\lambda$, where if $\lambda$ cuts successively adjacent regions $R'_i, i = 1, \ldots, k+1$ then $\alpha_i$ is the label of the common side of $R'_i$ and $R'_{i+1}$ which is *inside* $R'_{i+1}$. The sequence thus obtained is called the $\tilde{G}$–*cutting sequence* of $\lambda$. With the above labelling conventions, if $h \in \Gamma$ and $z \in \Delta$, then one can verify that the $\tilde{G}$–cutting sequence of the oriented geodesic from $z$ to $h(z)$ is a word in the generators $\tilde{G}_0$ representing $h$, see [1] for details.

We define $G$– and $\Gamma$–cutting sequences similarly. It is clear that the $G$–*cutting sequence* of the projection of the segment $\lambda$ to $\Omega$ is obtained from the $\tilde{G}$–sequence by omitting the labels $B^\pm$ and replacing $X'$ by $X$ and $Y'$ by $Y$. This specifies implicitly that the map $h \colon \pi_1(S) \to G = \pi_1(\mathbf{H}^3/G_\rho)$ is $h(X') = X, h(Y') = Y_\rho, h(B') = \text{id}$. Likewise the $\Gamma$–cutting sequence of the projection of $\lambda$ to $\mathbf{C} - \mathcal{L}$ is obtained from the $\tilde{G}$–sequence by replacing the labels $B'^\pm$ with $\beta'^\pm$, $X'$ by $\xi$ and $Y'$ by $\eta$. Clearly, since the combinatorics of all three diagrams in figure 2 are the same, we can read off the $G$–sequence from the $\Gamma$–sequence by omitting the labels $\beta^\pm$ and replacing $\xi$ by $X$ and $\eta$ by $Y$. This is a key point in our idea.

In practice, the cutting sequence is read off from weighted diagram by a simple combinatorial procedure. For definiteness, suppose we have a weighted diagram on the region $R''$ as in figure 3 of [4]. First, redraw the diagram replacing the line $l(\alpha, \beta)$ with weight $n(\alpha, \beta)$ by $n(\alpha, \beta)$ parallel and pairwise disjoint arcs joining the sides $\alpha, \beta$. As explained above, these arcs link in a unique order to form a simple closed loop $\lambda$ on $S$. Pick an orientation and initial point on $\lambda$. To follow the convention described above, every time $\lambda$ crosses an edge $s$ of $R''$, write down the label on $s$ and *outside* $R''$. Thus, if an oriented arc of $\lambda$ has initial point on an edge labelled $\alpha$ inside $R''$ and final point an edge labelled $\beta$ inside $R''$, then its contribution to the cutting sequence is $\alpha, \beta^{-1}$. The cutting sequence thus obtained is a word in the generators of the group $\Gamma$. Changing the initial point of $\lambda$ cyclically permutes the cutting sequence, so that the corresponding words are conjugate elements in $\Gamma$, while reversing the orientation of $\lambda$ produces the inverse word. Thus the loop $\lambda$ defines an I–equivalence class in $\Gamma$.

Now let $p/q \in \mathbf{Q} \cup \infty$ and let $L_{p/q}$ denote some line of slope $p/q$ in $\mathbf{C}$. Its $\Gamma$–cutting sequence is periodic, and the word $V_{p/q} \in G$ of proposition 2.2 representing $h(\gamma(p/q)) \in G$ is obtained by the procedure described above. Notice that $V_{p/q}$ is automatically cyclically reduced. Clearly, $h(\gamma(p/q))$ is equally represented by the word $V_{p/q}^{-1}$ corresponding to the cutting sequence of the line $L_{p/q}$ with its orientation reversed.

The remark on page 77 of [4] gives some examples. We note that the words given in the text are correct, but should be read off relative to the corrected labelling of figure 2 in which $Y$ and $Y^{-1}$ are interchanged.

That the words $V_{p/q}$ are defined only up to cyclic conjugation and inversion is another source of confusion in [4]. Only the I–equivalence class is well defined. As noted above, this equivalence class should also not change when $L_{p/q}$ is replaced by a parallel line of the same slope. In fact, it is clear that there are only a finite set of possible cutting sequences obtained from parallel translates of a line segment of finite length and that these sequences differ only by cyclic permutation. We denote the I–equivalence class in $G$ thus obtained by $C_{p/q}$.

## 1.1  The enumeration error

In accordance with the comments in the introduction, our task is to identify when two equivalence classes $C_{p/q}$ and $C_{r/s}$ coincide. This problem is discussed in remark 2.5 on page 79 of [4], where it is stated correctly that $V_{p/q} \sim V_{r/s}$ if $r/s = p/q + 2n, n \in \mathbf{Z}$. However, the claim in that remark that if $0 \le p/q <$ $r/s < 2$ then $\gamma(p/q)$ and $\gamma(r/s)$ are distinct is wrong; in fact, as explained in the proof of theorem 1.2 below, only $q > 0$ and $|p|$ are invariants of the class $C_{p/q}$. Thus, contrary to the claims implicit in [4], we have:

**Lemma 1.1** For $p/q \in \mathbf{Q}$, the classes $C_{p/q}$ and $C_{-p/q}$ coincide.

**Proof** Let $L_{p/q}$ be a line of rational slope $p/q \in \mathbf{Q}$ with initial point on the edge of $R''$ joining vertices $0, i$. Its reflection $L_{-p/q}$ in the imaginary axis has slope $-p/q$; let $V_{\pm p/q}$ be the words obtained from the $G$-cutting sequences of $L_{\pm p/q}$ as above. It is easy to see that the $\Gamma$-sequences of $L_{\pm p/q}$ differ by interchanging $\xi$ with $\xi^{-1}$, $\eta$ with $\eta^{-1}$, and $\beta$ with $\beta^{-1}$. (The interchange of $\beta$ with $\beta^{-1}$ happens because in the tesselation of $\mathbf{C} - \mathcal{L}$ by images of $R''$ under $\Gamma$, the labels $\beta$ and $\beta^{-1}$ alternate along horizontal lines.) Therefore $C_{p/q}(X^{-1}, Y^{-1}) = C_{-p/q}(X, Y)$.

Now compare two lines of the same slope $p/q$ which differ by vertical translation by $i$. Their cutting sequences differ by interchanging $\xi$ with $\xi^{-1}$, $\eta$ with $\eta^{-1}$, and $\beta$ with $\beta^{-1}$; in addition the position of the $\beta$ terms in the sequence shifts relative to that of the $\xi's$ and $\eta's$. (For example the sequence for $1/1$ with initial point between $0$ and $i$ is $\xi\eta^{-1}\beta^{-1}$, while with initial point between $i$ and $2i$ we get $\xi^{-1}\beta\eta$.) Since the position of $\beta^\pm$ relative to $\xi^\pm, \eta^\pm$ does not affect the $X, Y$ sequence, we get $C_{p/q}(X^{-1}, Y^{-1}) = C_{p/q}(X, Y)$.

Combining these observations gives the proof. □

We also need to know there are no other identifications. We have:

**Theorem 1.2** *The classes $C_{p/q}, C_{r/s}, p/q, r/s \in \mathbf{Q} \cup \infty$ are equivalent if and only if $r/s = p/q + 2n$ or $-r/s = p/q + 2n$, $n \in \mathbf{Z}$.*

**Proof** As discussed in remark 2.5 of [4] on page 79, a (left) Dehn twist about the curve $\gamma(\infty)$ represented by $\beta \in \Gamma$ induces an automorphism of $\Gamma$ which maps $\gamma(p/q)$ to $\gamma(2+p/q)$. Since this automorphism induces the identity on $G$, we have $C_{p/q} = C_{2+p/q}$. This can also be seen by representing $C_{p/q}$ and $C_{2+p/q}$ by the cutting sequences of lines of slope $p/q$ and $2 + p/q$ in $\mathbf{C}$. Reading off the two cutting sequences starting from the same inital point, it is easy to see that, while the $\Gamma$–sequences differ, the induced $G$–sequences are the same.

To complete the proof, it only remains to show that if $r/s \neq \pm(p/q + 2n)$ then $C_{p/q}, C_{r/s}$ are distinct.

As stated in remark 2.5 of [4], the cutting sequence of $L_{p/q}$ has length $2q$. Moreover, the words $V_{p/q}$ and $V_{r/s}$ are cyclically reduced and so, since $G$ is a free group, are conjugate only if they have the same length. Thus, since we are assuming that $q, s \geq 0$, a necessary condition for $C_{p/q} = C_{r/s}$ is that $q = s$.

It is also stated in remark 2.5 that $p$ can be deduced from number of sign changes in the exponents of $X$ and $Y$ in $V_{p/q}$. This is not quite correct, and herein lies the root of error number 1. Since the number of sign changes is necessarily non-negative, one verifies that in all cases we can only obtain $|p|$ and not $p$ from $C_{p/q}$, in other words, $C_{p/q} = C_{r/s}$ implies $|p| = |r|$ but, contrary to the claim of remark 2.5, the number of sign changes cannot be used to distinguish the classes of $V_{p/q}$ and $V_{-p/q}$. This is correct, since by lemma 1.1 the two classes $C_{p/q}$ and $C_{-p/q}$ coincide. □

**Remark 1.3** For future reference, we note that a similar argument to the above shows that $V_{1+p/q}$ can be obtained from $V_{p/q}$ by interchanging $Y$ and $Y^{-1}$, more precisely, that $V_{1+p/q}(X,Y) = V_{p/q}(X,Y^{-1})$.

This completes the discussion of the first error.

## 2 Connectivity

Let $g \in G_\rho$ correspond to a simple closed geodesic $\gamma$ on $\partial \mathcal{C}/G$. The trace $\operatorname{Tr} g$ is a polynomial in $\rho$ with integer coefficients. It is claimed in [4] theorem 4.1 that the pleating ray $\mathcal{P}_\gamma$ has a unique connected component with a unique endpoint $\rho^* = \rho^*(\gamma)$ on $\partial \mathcal{R}$. At this endpoint, $\operatorname{Tr} g = \operatorname{Tr} g(\rho^*) = -2$ ( [4] proposition 4.2) and $g$ is an accidental parabolic. The group $G_\rho$ is free, $\Omega(G_\rho) \neq \emptyset$, and therefore $G_\rho$ is maximally parabolic as in [6], i.e., $G_{\rho^*}$ contains the maximal number of rank 1 parabolic subgroups among subgroups of $PSL(2, \mathbf{C})$ isomorphic to $G_{\rho^*}$.

The map $\rho \mapsto \bar\rho$ induces the maps $X \mapsto X, Y_\rho \mapsto Y_{\bar\rho}$ and hence an isomorphism $J: G_\rho \to G_{\bar\rho}$; clearly, $J$ is type preserving, i.e., $g \in G_\rho$ and $J(g) \in G_{\bar\rho}$ are either both parabolic or both loxodromic It is also clear that $\rho \mapsto \bar\rho$ maps $\mathcal{R}$ to itself. We note that this does not contradict the uniqueness of maximally pinched groups asserted in theorem III of [6], because the conjugation $G_{\rho^*} \to G_{\bar\rho^*}$ is *anti*holomorphic. However, it does mean that $\bar\rho^*$ should be also be an endpoint of the pleating ray $\mathcal{P}_\gamma$, which contradicts theorem 4.1 of [4].

This contradiction is resolved by the corollary to the following lemma.

**Lemma 2.1** *Let $j\colon \hat{\mathbf{C}} \to \hat{\mathbf{C}}$ be an conformal or anticonformal bijection and let $j_*(M) = jMj^{-1}, M \in SL(2, \mathbf{C})$. Let $\langle G_\rho; X, Y_\rho \rangle$ be the marked free group with generators $X, Y_\rho$. Suppose that the pleating locus $pl(G_\rho)$ consists of a simple closed geodesic represented by a word $W(X, Y_\rho)$. Let $j_*(G_\rho)$ denote the marked group with generators $j_*(X), j_*(Y_\rho)$. Then $pl(j_*(G_\rho))$ is a closed geodesic represented by the word $W(j_*(X), j_*(Y_\rho))$.*

**Proof** Recall from [4] that a subgroup of a Kleinian group $G$ is called *F–peripheral* if it is Fuchsian and if one of the two open round discs bounded by its limit set contain no points of $\Lambda(G)$. If a geodesic $\gamma$ is the pleating locus of $G_\rho, \rho \in \mathcal{R}$, then, as discussed on page 82 of [4], $\gamma$ divides $\partial \mathcal{C}/G_\rho$ into two connected components each of which is a sphere with two punctures and one hole. The lifts of these components to $\mathbf{H}^3$ lie in hyperbolic planes which separate $\mathbf{H}^3$ into two open hyperbolic half spaces, one of which contains no points of $\mathcal{C}$. The half spaces meet $\hat{\mathbf{C}}$ in open discs which have empty intersection with the limit set $\Lambda(G_\rho)$, so that the subgroups of $G_\rho$ which leave these discs invariant are F–peripheral and contain the element suitable conjugates of $W(X, Y_\rho)$. In particular, $W(X, Y_\rho)$ lies in in two non-conjugate F–peripheral subgroups of $G_\rho$.

We showed in [4] proposition 3.6 that conversely, if an element $g \in G_\rho$ which represents a simple closed geodesic $\gamma$ on $\partial \mathcal{C}/G$ lies in two non-conjugate F–peripheral subgroups, then $pl(G_\rho) = \gamma$.

Now since the map $j$ is conformal or anticonformal, it maps circles to circles. It is also clear that $\Lambda(j_*(G_\rho)) = j(\Lambda(G_\rho))$, and that if a subgroup of $G_\rho$ is Fuchsian, so is its image in $j_*(G_\rho)$. Hence $j_*$ preserves F–peripheral subgroups.

Thus $W(X, Y_\rho)$ lies in two non-conjugate F–peripheral subgroups of $G_\rho$ if and only if $j_*(W(X, Y_\rho)) = W(j_*(X), j_*(Y_\rho))$ lies in two non-conjugate F–peripheral subgroups of $j_*(G_\rho)$. The result follows. $\square$

**Corollary 2.2** *Let $p/q \in \mathbf{Q} \cup \infty$ and let $j\colon \rho \mapsto \bar\rho$ be complex conjugation. Then $\mathcal{P}_{p/q} = j(\mathcal{P}_{p/q})$.*

**Proof** Let $\rho \in \mathcal{P}_{p/q}$ so that $pl(G_\rho) = \gamma_{p/q}$. As above, in the marked group $\langle G_\rho; X, Y_\rho \rangle$, the class $\gamma_{p/q}$ is represented by the word $V_{p/q}(X, Y_\rho)$ in $X, Y_\rho$. We apply lemma 2.1 to $j$ and compute $j_*(X) = X$ and $j_*(Y_\rho) = Y_{\bar\rho}$. Thus $j_*(G_\rho)$ is the marked group $\langle G_{\bar\rho}; X, Y_{\bar\rho} \rangle$, and $pl(G_{\bar\rho})$ is represented by the word $V_{p/q}(X, Y_{\bar\rho})$ which corresponds in the marked group $G_{\bar\rho}$ to $\gamma_{p/q}$. $\square$

**Remark 2.3** We can also apply lemma 2.1 to the involution $k\colon \rho \mapsto -\rho$. We find $k_*(X) = X^{-1}$ and $k_*(Y_\rho) = Y_{-\rho} = Y_\rho^{-1}$. Thus $k_*(G_\rho)$ is the marked group $\langle G_{-\rho}; X^{-1}, Y_\rho^{-1}\rangle$. If as above, $pl(G_\rho)$ is represented by the word $V_{p/q}(X, Y_\rho)$ in $X, Y_\rho$, then $pl(G_{-\rho})$ is represented by the word $V_{p/q}(X^{-1}, Y_\rho^{-1})$. Now as in lemma 1.1 and remark 1.3 above, $V_{p/q}(X^{-1}, Y_\rho^{-1}) \sim V_{p/q}(X, Y_\rho)$ and $V_{p/q}(X, Y_\rho) \sim V_{1+p/q}(X, Y_\rho^{-1})$. Thus the pleating locus of the *marked* group $\langle G_{-\rho}; X, Y_{-\rho}\rangle = \langle G_{-\rho}; X, Y_\rho^{-1}\rangle$ is $\gamma_{1+p/q}$ which, using remark 2.3, is the same as $\gamma_{1-p/q}$.

Although as groups $G_\rho$ and $k(G_\rho)$ are the same, $k_*$ is *not* the the standard isomorphism and $k_*(V_{p/q}(X, Y_\rho)) \ne V_{p/q}(X, Y_{-\rho})$. This explains why the endpoints of the rays $\mathcal{P}_{p/q}$ and $k(\mathcal{P}_{p/q}) = -\mathcal{P}_{p/q}$ correspond to *different* maximally pinched groups.

## 2.1 The connectivity error

We can now prove a correct form of theorem 4.1 of [4]. Recall that the hyperbolic locus of the trace polynomial $\operatorname{Tr} V_{p/q}$ is the set

$$\tilde{\mathcal{H}}_{p/q} = \{\rho \in \mathbf{C} : \Im \operatorname{Tr} V_{p/q} = 0, |\Re \operatorname{Tr} V_{p/q}| > 2\},$$

and that $\mathcal{P}_{p/q}$ is a union of connected components of $\tilde{\mathcal{H}}_{p/q}$.

**Theorem 2.4** For $0 < p/q < 1$, the rational pleating ray $\mathcal{P}_{p/q}$ consists of exactly two connected components of the hyperbolic locus $\tilde{\mathcal{H}}_{p/q}$. These rays are the branches which asympotically have arguments $-e^{\pi ip/q}$ and $-e^{-\pi ip/q}$. They are complex conjugate 1–manifolds, with unique and complex conjugate endpoints on $\partial \mathcal{R}$.

**Proof**  In [4], we argued by "induction on the Farey tree". Once again, there is an error in the argument which can be corrected using corollary 2.2.

For the rays $\mathcal{P}_{0/1}$ and $\mathcal{P}_{1/1}$ we argue exactly as in [4] proposition 3.8. (The assertion that these special rays have one connected component is correct; we note that since they are contained the real axis, they are invariant under complex conjugation so the contradiction explained above does not occur.)

Now suppose we have the result for $\mathcal{P}_{p/q}$ and $\mathcal{P}_{r/s}$ for which $ps - rq = -1$. Let $\mathbf{H}^+$ and $\mathbf{H}^-$ denote the upper and lower half planes respectively. The argument in [4] shows that an arc in $\mathbf{H}^+$ joining the components $\mathcal{P}_{p/q}^+$ to $\mathcal{P}_{r/s}^+$ of $\mathcal{P}_{p/q}$ to $\mathcal{P}_{r/s}$ in $\mathbf{H}^+$ must intersect $\mathcal{P}_{(p+r)/(q+s)}$. Also as in [4], the only branch of $\tilde{\mathcal{H}}_{(p+r)/(q+s)}$ whose asymptotic direction lies between directions $-e^{\pi i p/q}$ and $-e^{\pi i r/s}$ is the one with asymptotic direction $-e^{\pi i(p+r)/(q+s)}$, and this must therefore be coincident with a component of $\mathcal{P}_{(p+r)/(q+s)}$. Similarly an arc in $\mathbf{H}^-$ joining the components $\mathcal{P}^-_{p/q}$ to $\mathcal{P}^-_{r/s}$ of $\mathcal{P}_{p/q}$ to $\mathcal{P}_{r/s}$ in $\mathbf{H}^-$ must intersect a component of $\mathcal{P}_{(p+r)/(q+s)}$, with asymptotic direction $-e^{-\pi i(p+r)/(q+s)}$. This gives the result.

(The problem with the argument in [4] is that we forgot to consider arcs joining the components of $\mathcal{P}^+_{p/q}$ to $\mathcal{P}^+_{r/s}$ in $\mathbf{H}^+$ and running through $\mathbf{H}^-$.)

Notice that the picture obtained in this way is entirely consistent with remark 2.3 above. □

## 3  Conclusion

The errors above do not substantially effect any of the conclusions of [4]. Theorems 1.1 and 2.4 have obvious extensions to irrational laminations, which we shall not spell out here. The only other result which is changed in consequence of the errors is theorem 5.4.

In [4], to deal with irrational rays $\lambda \in \mathbf{R}$, we introduced the complex pleating length $L_\lambda$, and referred to the methods of [3], section 7.1 to show that these rays were 1–manifolds with the connectivity claimed. In fact, the argument in [3] has a gap: we omitted to show that the pleating variety $\mathcal{P}_\lambda$ is open in the real locus of $L_\lambda$. This crucial fact is proved in a more general context in [5]. For a corrected version of the arguments required in a one dimensional parameter space, we refer to [8]. We note also that by the improved techniques of [5], it follows that even on irrational rays $\lambda \in \mathbf{R}$, the range of the complex pleating length $L_\lambda$ (see [4] page 88) is $\mathbf{R}^+$.

Let $j$ denote complex conjugation and define an equivalence relation on $\mathbf{R}$ by $x \approx y$ if and only if $x = \pm y + 2n, n \in \mathbf{Z}$. We can think of the pleating locus $pl(\rho)$ as a $\approx$–equivalence class in $\mathbf{R}$. Then the map

$$\Pi \colon \mathcal{R} \to \mathbf{R}/\approx \times \mathbf{R}^+, \ \Pi(\rho) = (pl(\rho), L_{pl(\rho)}(\rho)),$$

factors through $j$. We denote the induced map, $\tilde{\Pi}$.

We obtain:

**Theorem 3.1**  *The map*

$$\tilde{\Pi} \colon \mathcal{R}/j \to \mathbf{R}/\approx \times \mathbf{R}^+$$

*is a homeomorphism.*

*Department of Mathematics, Osaka City University, Osaka 558, Japan*
and
*Mathematics Institute, Warwick University, Coventry CV4 7AL, England*

Email: `komori@sci.osaka-cu.ac.jp, cms@math.warwick.ac.uk`